\documentclass[12pt]{article}
\usepackage{amsmath}
\usepackage{amsbsy}
\usepackage{amssymb}
\usepackage[dvips]{graphicx}
\newtheorem{theo}{Theorem}[section]

\newcommand{\s}{{\mbox{\boldmath $s$}}}

\def\re{\mathbb R}

\def\h{\mathbb H}
\def\s{\mathbb S}

\def\lra{\longrightarrow}

\def\qed{\hfill $\sqcap \hskip-6.5pt \sqcup$}
\DeclareMathOperator{\diver}{\text{div}}

\title{A Halfspace Theorem for
    Mean Curvature  \\
$H=\frac{1}{2}$ surfaces
in ${\mathbb H}^2\times {\mathbb R}$}

\author{Barbara Nelli - Ricardo Sa Earp}

\date{}

\begin{document}

\maketitle

\begin{abstract}

 {\em We prove a vertical halfspace
theorem for surfaces with constant mean curvature $H=\frac{1}{2},$
properly immersed  in the product space $\h^2\times\re,$ where
$\h^2$ is the hyperbolic plane and $\re$ is the set of real
numbers. The proof is  a geometric application of the classical
maximum principle for second order elliptic PDE, using the family
of non compact rotational $H=\frac{1}{2}$ surfaces in
$\h^2\times\re.$ }
\end{abstract}


\section{ Introduction}
\label{introduction}

This is a revised version of the article that we submit before.
There was  a problem in the construction of graphical ends. We are
presently working to fix it (replace the previous boundary with a
planar boundary curve and use  Perron method). The main  geometric
constructions will be mantained. Here we present the halfspace
type theorem, that correspond to Section 4 of the previous
article.

D. Hofmann e W. Meeks proved a beautiful theorem on minimal
surfaces, the so-called   "Halfspace Theorem" in \cite{HM}: there
is no non planar, complete,  minimal surface properly immersed in
a halfspace of $\re^3.$  We focus in this paper complete surfaces
with constant mean curvature $H=\frac{1}{2}$ in the product space
$\h^2\times\re,$ where $\h^2$ is the hyperbolic plane and $\re$ is
the set of real numbers. In the context of $H$-surfaces in
$\h^2\times\re,$ it is natural to investigate about halfspace type
results.

Before stating our result we would like to emphasize that,
in last years there has been work in  $H$-surfaces in homogeneous
3-manifolds, in particular in the product space $\h^2\times\re:$ new examples were produced and many
theoretical results as well.

Halfspace theorem for minimal surfaces in $\h^2\times\re$ is
false, in fact there are many  vertically bounded complete minimal
surfaces in $\h^2\times\re$ \cite{ST2}.  On the contrary, we are
able to prove the following result for $H=\frac{1}{2}$ surfaces.

\begin{theo} \label{halfspacetheorem}
 Let $S$ be a simply connected rotational
surface with constant mean curvature $H=\frac{1}{2}.$  Let $\Sigma$
be a complete surface  with constant mean curvature $H=\frac{1}{2},$
different from a rotational simply connected one. Then, $\Sigma$ can
not be properly immersed in the mean convex side of $S.$
\end{theo}

In \cite{HRS1}, L. Hauswirth, H. Rosenberg and J. Spruck  prove a
halfspace type theorem for surfaces on one side of a horocylinder.

 The result in \cite{HRS1} is different in nature from our result because in \cite{HRS1}, the
 "halfspace" is one side of a horocylinder, while for us, the "halfspace" is the mean convex  side of the
 rotational simply connected surface.

The proof of our result  is  a geometric  application of the classical maximum principle  to
surfaces  with constant mean curvature $H=\frac{1}{2}$ in $\h^2\times\re.$

\vspace{.5cm}

{\bf Maximum Principle.} {\em Let $S_1$ and $S_2$ be two connected surfaces of constant mean curvature $H=\frac{1}{2}.$ Let
$p\in S_1\cap S_2$ be a point such that $S_1$ and $S_2$ are tangent at $p,$ the mean curvature vectors of $S_1$ and $S_2$ at $p$ point towards the same side and $S_1$ is on one side of $S_2$  in a neighborhood  of $p.$ Then $S_1$ coincide with $S_2$ around  $p.$  By analytic continuation, they coincide everywhere.}
\vspace{.5cm}

The proof of the Maximum Principle is based on the fact that a constant mean curvature surface in $\h^2\times\re$ locally satisfies a second order elliptic PDE (cf. \cite{H}, \cite{GT}  where the author prove the Maximum Principle in $\re^n;$ the proof generalizes to space forms and to $\h^2\times\re$ as well).

We notice that our surfaces are not compact, while the classical maximum principle applies at a finite point.
It will be clear in the proof of  Theorem \ref{halfspacetheorem} that  we that  we are able to reduce the
analysis to finite tangent points,
because of  the geometry of rotational surfaces of constant mean curvature $H=\frac{1}{2}.$

Our halfspace Theorem leads to the following conjecture (strong halfspace theorem).
\vspace{.5cm}

{\bf Conjecture.} {\em Let $\Sigma_1,$ $\Sigma_2$
be two complete properly embedded surfaces with constant mean curvature
$H=\frac{1}{2},$ different from the rotational simply connected one.
Then $\Sigma_i$ can not lie in the mean convex side of $\Sigma_j,$
$i\not=j.$}
\vspace{.5cm}

For $H>\frac{1}{\sqrt 2}$ the conjecture is true and it is known as maximum principle at infinity (cf. \cite{NR1}).

\section{Vertical Halfspace Theorem}
\label{halfspacesection}

R. Sa Earp and E. Toubiana find  explicit integral formulas for
rotational surfaces of constant mean curvature $H\in(0,\frac{1}{2}]$
in \cite{ST}.   A careful description of the geometry of these
surfaces is contained in the Appendix of \cite{NSST}.

Here we recall some properties of  rotational surfaces of constant mean curvature $H=\frac{1}{2}.$
For any $\alpha\in\re_+,$ there exists a rotational surface ${\cal H}_{\alpha}$ of constant mean curvature $H=\frac{1}{2}.$

For $\alpha\not=1,$  the surface ${\cal H}_{\alpha}$ has two vertical ends (where a vertical
end  is   a topological  annulus,  with no asymptotic point at
finite height) that are graphs over the exterior of a disk $D_{\alpha}$ of hyperbolic radius $r_{\alpha}=|\ln\alpha|.$

By graph we mean the following:   the graph of a function $u$ defined on a subset $\Omega$ of $\h^2$ is
$\left\{(x, y, t)\in\Omega\times\re \ \vert \ t=u(x,y)\right\}.$
When the graph has constant mean curvature $H,$
$u$   satisfies the following second order elliptic PDE

\begin{equation}
\label{graph} \diver_{\h}\left(\frac{\nabla_{\h} u}{W_u}\right)=2H
\end{equation}

where $\diver_{\h},$ $\nabla_{\h} $ are the hyperbolic divergence
and gradient respectively and $W_u=\sqrt{1+|\nabla_{\h}u|_{\h}^2},$
being  $|\cdot|_{\h}$ the norm in $\h^2\times\{0\}.$

Furthermore, up to vertical translation, one can assume that ${\cal H}_{\alpha}$ is symmetric with respect to the  horizontal plane $t=0.$
For $\alpha=1,$ the surface ${\cal H}_1$   has only one end and it is a graph over $\h^2$ and it is denoted by $S.$

When $\alpha>1$ the surface  ${\cal H}_{\alpha}$ is not embedded.
The self intersection set is a horizontal circle
on the plane $t=0.$ For  $\alpha<1$ the surface  ${\cal H}_{\alpha}$ is  embedded.

 For any $\alpha\in\re_+,$ each  end of the   surface  ${\cal H}_{\alpha}$ is the vertical graph of a
function $u_{\alpha}$  over  the exterior of a disk $D_{\alpha}$ of radius $r_{\alpha}.$  The asymptotic
behavior of $u_{\alpha}$ has the following form:  $ u_{\alpha}(\rho)\simeq \frac{1}{\sqrt{\alpha}} e^{\frac{\rho}{2}},$ $\rho\lra\infty,$ where $\rho$ is the hyperbolic distance from the origin. The positive number $\frac{1}{\sqrt\alpha}\in\re_+$ is called the {\em growth} of the end.

\begin{figure}
\label{profilecurvefig}
\begin{center}
\includegraphics[scale=0.30]{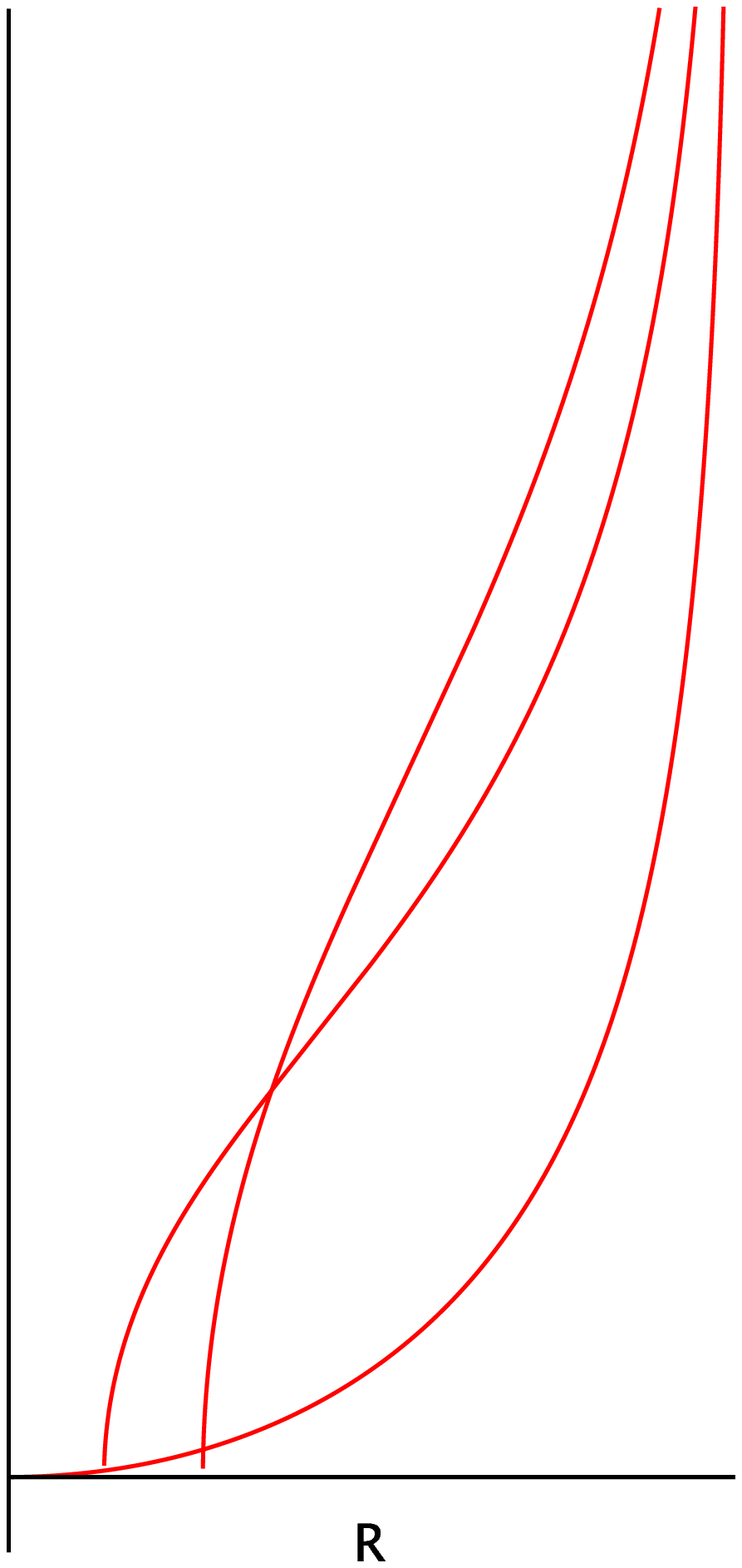}
\includegraphics[scale=0.50]{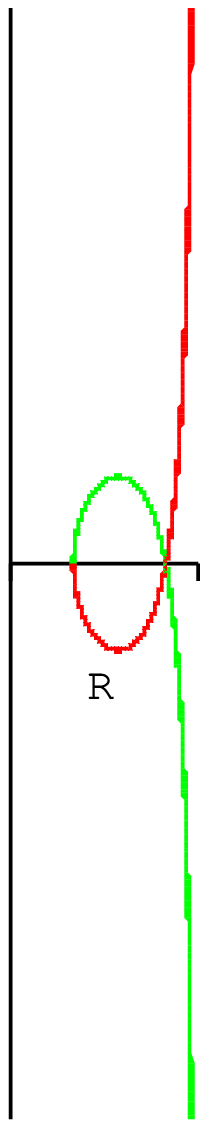}
 \caption{$H=\frac{1}{2}:$ the profile curve in the embedded and
immersed case ($R=\tanh\rho$).}
\end{center}
\end{figure}

The function $u_{\alpha}$
is  vertical along the  boundary of $D_{\alpha}.$ Furthermore the radius
$r_{\alpha}$   is always greater or equal to zero, it is zero if and
only if $\alpha=1$ and tends to infinity as $\alpha\lra 0$ or $\alpha\lra \infty.$ As we
pointed out before,  the function $u_{1}=2\cosh \left(\frac{\rho}{2}\right)$
is entire and its graph
corresponds to the unique simply connected example $S.$

Notice that, any end of an immersed rotational surface ($\alpha>1$) has growth  smaller than the growth of  $S,$ while any end of an embedded rotational surface ($\alpha<1$) has  growth greater than the growth of $S.$

Theorem \ref{halfspacetheorem} is called "vertical" because the end of the surface $\Sigma$ is vertical, as it is contained in the mean convex side of $S.$

\vspace{.5cm}

 {\bf Proof of Theorem \ref{halfspacetheorem}.} One can assume that the surface $S$ is tangent to
the
 slice $t=0$ at the origin and it is contained in $\{t\geq
 0\}.$
 Suppose, by contradiction, that $\Sigma$ is contained in the mean
 convex side of $S.$
 Lift  vertically $S.$  If there is an interior contact point
 between $\Sigma$ and the translation of $S,$ one has a
 contradiction  by the maximum principle. As $\Sigma$ is properly immersed,
 $\Sigma$ is asymptotic at infinity  to a vertical translation of $S.$
 One can assume that
 the surface $\Sigma$ is asymptotic to the $S$  tangent to the
 slice $t=0$ at the origin and  contained in $\{t\geq
 0\}.$

Let $h$ be the height of  one lowest point of $\Sigma.$ Denote by
$S(h)$ the vertical lifting of $S$ of ratio $h.$ One has one of the
following facts.
\begin{itemize}
\item $S(h)$ and $\Sigma$ has a first finite contact point $p:$  this means that
$S(h-\varepsilon)$ does not meet $\Sigma$ at a finite point, for
$\varepsilon>0$ and then $S(h)$ and $\Sigma$ are tangent   at $p$
with mean curvature vector pointing in the same direction. In this
case, by the maximum principle $S(h)$ and $\Sigma$ should coincide.
Contradiction.

\item $S(h)$ and $\Sigma$ meet  at  a point $p,$
but $p$ is not a first contact point.Then,  for $\varepsilon$ small
enough, $S(h-\varepsilon)$ intersect $\Sigma$ transversally.
\end{itemize}

  Denote by  $W$  the non compact
subset of $\h^2\times\re$ above $S$ and below $S(h-\varepsilon).$

It follows from the maximum principle that there are no compact
component of $\Sigma$ contained in $W.$ Denote by $\Sigma_1$ a non
compact connected component of $\Sigma$ contained in $W.$ Note that the
boundary of $\Sigma_1$ is contained in $S(h-\varepsilon).$ Consider
the family of  rotational  non embedded surfaces ${\cal
H}_{\alpha},$ $\alpha>1.$ Translate each ${\cal H}_{\alpha}$
vertically in order to have the waist on the plane
$t=h-\varepsilon.$ By abuse of notation, we continue to call the
translation, ${\cal H}_{\alpha}.$ The surface ${\cal H}_{\alpha}$
intersects the plane $t=h-\varepsilon$ in two circles. Denote by
$\rho_{\alpha}$ the radius of the larger circle. Denote by ${\cal
H}^+_{\alpha},$ the part of the surface outside the cylinder of
radius $\rho_{\alpha}.$ Notice that ${\cal H}^+_{\alpha}$ is
embedded.  By the geometry of the ${\cal H}^+_{\alpha},$ when
$\alpha$ is great enough, say $\alpha_0, $ ${\cal H}^+_{\alpha_0}$
is outside the mean convex side of $S.$  Then,
${\cal H}^+_{\alpha_0}$ does not intersect $\Sigma.$ Furthermore,
when $\alpha\longrightarrow 1,$ ${\cal H}^+_{\alpha}$ converge to
$S(h-\varepsilon).$  Now, start to decrease $\alpha$ from $\alpha_0$
to one. Before reaching $\alpha=1,$ the surface ${\cal
H}^+_{\alpha}$ first meets $S$ and then  touches $\Sigma_1$ tangentially at an
interior  finite point, with $\Sigma_1$  above ${\cal H}^+_{\alpha}.$  This depends on the following two facts.

\begin{itemize}
\item  The boundary of $\Sigma_1$ lies on $S(h-\varepsilon)$ and the
boundary of any of the ${\cal H}^+_{\alpha}$ lies on the horizontal
plane $t=h-\varepsilon.$

\item The growth of any of the ${\cal H}^+_{\alpha}$ is strictly smaller than
the growth of $S.$ Thus the end of ${\cal H}^+_\alpha$ is outside
the end of $S$.

\end{itemize}

The existence of an such interior tangency  point is a contradiction by the
maximum principle.

\qed

\textsc{Barbara Nelli}

{\em  Universit\'a di L'Aquila

nelli@univaq.it}

\textsc{Ricardo Sa Earp}

{\em PUC, Rio de Janeiro

earp@mat.puc-rio.br}

\end{document}